\documentclass[graybox]{svmult}

\usepackage{mathptmx}       % selects Times Roman as basic font
\usepackage{helvet}         % selects Helvetica as sans-serif font
\usepackage{courier}        % selects Courier as typewriter font
\usepackage{type1cm}        % activate if the above 3 fonts are
\usepackage{makeidx}         % allows index generation
\usepackage{graphicx}        % standard LaTeX graphics tool
\usepackage{multicol}        % used for the two-column index
\usepackage[bottom]{footmisc}% places footnotes at page bottom

\makeindex             % used for the subject index
\usepackage{amssymb,amsmath,mathrsfs,amsfonts,graphicx,epsf}
\usepackage{color}
\usepackage{xcolor}
\usepackage[active]{srcltx}

\usepackage{hyperref}            % liens en pdf avec pdfLatex (pour les Mac)

 \usepackage{esint}

 \usepackage{fullpage}

\newcommand{\ord}{\mathrm{ord}}
\newcommand{\qon}{\overline{Q}_N}

\newcommand{\bD}{{\cal D}}
\renewcommand{\wp}{\hat \omega^p}
\newcommand{\mup}{\hat \mu^p}
\newcommand{\muq}{\hat \mu^q}

\renewcommand{\span}{\mathrm{span}}

\newcommand{\bt}{\begin{theorem}}
\newcommand{\et}{\end{theorem}}
\newcommand{\bl}{\begin{lemma}}
\newcommand{\el}{\end{lemma}}
\newcommand{\bp}{\begin{proposition}}
\newcommand{\ep}{\end{proposition}}
\newcommand{\bc}{\begin{corollary}}
\newcommand{\ec}{\end{corollary}}
\newcommand{\bdeff}{\begin{definition}}
\newcommand{\edeff}{\end{definition}}
\newcommand{\brem}{\begin{remark}}
\newcommand{\erem}{\end{remark}}
\renewcommand{\r}[1]{(\ref{#1})}
\newcommand{\con}{{\mathcal C}}

\def\diam{\mathop{\mathrm{diam}}}

\def\ss{\mathcal{S}}
\newcommand{\hh}{{\mathcal H}}

\newcommand{\qn}{{Q_N}}

\newcommand{\bi}{\begin{itemize}}
\newcommand{\iii}{\item}
\newcommand{\ei}{\end{itemize}}
\newcommand{\bd}{\begin{description}}
\newcommand{\ed}{\end{description}}

\newcommand{\bqn}{\begin{eqnarray}}
\newcommand{\eqn}{\end{eqnarray}}
\newcommand{\eqnn}{\nonumber\end{eqnarray}}

\newcommand{\ba}[1]{\begin{array}{#1}}
\newcommand{\ea}{\end{array}}

\newcommand{\R}{\mathbb{R}}

\newcommand{\N}{\mathbb{N}}

\newcommand{\eps}{\epsilon}

\newcommand{\VecM}{\mathrm{Vec}(M)}

\begin{document}

\title*{Hausdorff measures and dimensions in non equiregular sub-Riemannian manifolds
}
\author{R.~Ghezzi and F.~Jean}
\institute{R.~Ghezzi \at Scuola Normale Superiore, Piazza dei Cavalieri 7, 56126 Pisa Italy, \email{roberta.ghezzi@sns.it}
\and F.~Jean \at ENSTA ParisTech, 828 Boulevard des Mar\'{e}chaux
91762 Palaiseau, France,  and
  Team  GECO, INRIA Saclay --
\^{I}le-de-France, \email{frederic.jean@ensta-paristech.fr}}
\maketitle
\abstract{This paper is a starting point towards computing the Hausdorff dimension of submanifolds and the Hausdorff volume of small balls in a sub-Riemannian manifold with singular points. We first consider the case of  a  strongly equiregular submanifold, i.e., a smooth submanifold $N$ for which the growth vector of the distribution $\bD$ and the growth vector of the intersection of $\bD$  with $TN$ are constant on $N$. In this case, we generalize the result in \cite{mitchell}, which relates the Hausdorff dimension to the growth vector of the distribution. We then consider analytic sub-Riemannian manifolds and, under the assumption that the singular point $p$ is typical, we state a theorem which characterizes the Hausdorff dimension of the manifold and the finiteness of the Hausdorff volume of small balls $B(p,\rho)$ in terms of the growth vector of both the distribution and the intersection of the distribution with the singular locus, and of the nonholonomic order at $p$ of the volume form on $M$ evaluated along some families of vector fields.}

\section{Introduction}

The main motivation of this paper arises from the study of sub-Riemannian manifolds as particular metric spaces. Recall that  a sub-Riemannian manifold is a triplet
 $(M,\bD,g)$, where  $M$ is a
smooth manifold, $\bD$ a Lie-bracket generating subbundle of
$TM$ and $g$ a   Riemannian metric on $\bD$.
 The absolutely continuous paths which are
almost everywhere tangent to $\bD$ are called horizontal and their  length  is
obtained as in Riemannian geometry integrating the norm of
their tangent vectors.  The sub-Riemannian
distance $d$ is defined as the infimum of length of horizontal paths between two given points.

 Hausdorff measures and spherical Hausdorff measures can be defined on sub-Riemannian manifolds using the sub-Riemannian distance.
 It is well-known that for these metric spaces the Hausdorff dimension is strictly greater than the topological one. Although the presence of an extra structure, i.e., the differential one, constitute a considerable help, computing Hausdorff measures and dimensions of sets is a difficult problem.
 In \cite{noi} we study Hausdorff measures of continuous curves, whereas in \cite{balu} the authors analyze the regularity of the top-dimensional Hausdorff measure   in the equiregular case (see the definition below). In the case of Carnot groups, Hausdorff measures of regular hypersurfaces have been studied in \cite{fssc} and in a more general context, a representation formula for the perimeter measure in terms of Hausdorff measure has been proved in \cite{ambro}.

In this paper we consider three questions: given a sub-Riemannian manifold $(M,\bD,g)$, $p\in M$ and a small $\rho>0$,
\bi
\iii[1.] what is the Hausdorff dimension $\dim_{H}(M)$?
\iii[2.] under which condition is the Hausdorff volume $\hh^{\dim_{H}(M)}(B(p,\rho))$  finite?
\iii[3.] the two preceding questions when $M$ is replaced by a submanifold $N$, i.e., what is $\dim_{H}(N)$ and when is $\hh^{\dim_{H}(N)}(N \cap B(p,\rho))$  finite?
\ei

A key feature to be taken into account is whether $p$ is regular or singular  for the sub-Riemannian manifold.
Given $i\geq 1$, define recursively  the
  submodule $\bD^i$ of $\VecM$ by
  $
\bD^{1}=\bD$, $\bD^{i+1}=\bD^{i}+[\bD,\bD^{i}]
$.
Denote by $\bD^i_p=\{X(p)\mid X\in\bD^i\}$. Since $\bD$ is Lie-bracket generating,
there exists $r(p)\in\N$ such that
$$
\{0\}=\bD^0_p\subset\bD_p^1\subset\dots\subset\bD_p^{r(p)}=T_pM.
$$
A point $p$ is regular  if, for every $i$, the dimensions $\dim\bD^i_q$ are constant as $q$ varies in a neighborhood of $p$. Otherwise, $p$ is said to be singular. A set  $S\subset M$ is equiregular if, for every $i$, $\dim\bD^i_q$ is constant as $q$ varies in $S$.
 For equiregular manifolds, questions 1 and 2 have been answered in  \cite{mitchell} (but with an incorrect proof, see \cite{montgomery} for a correct one). In that paper, the author shows that the Hausdorff dimension of an equiregular manifold $M$ is
 \begin{equation}\label{eqintro}
\dim_H(M)=Q, \quad \hbox{where } \sum_{i=1}^{r(p)}i(\dim\bD^i_p-\dim\bD^{i-1}_p),
 \end{equation}
and that the Hausdorff $Q$-dimensional measure near a regular point is absolutely continuous
with respect to any Lebesgue measure on $M$.
As a consequence, when $p$ is regular, the Hausdorff dimension of a small ball $B(p,\rho)$ is $Q$, and the Hausdorff $Q$-dimensional measure of $B(p,\rho)$ is finite.
%solving the first query.\red{reference for the second one?}
% and that  for every $p$
% \bi
%\iii[1.] what is the Hausdorff dimension $\dim_{H}(B(p,r))$?
%\iii[2.] under which condition is $\hh^{\dim_{H}}(B(p,r))$  finite?
%\ei

When there are singular points,
these problems have been mentioned in  \cite[Section~1.3.A]{gromovcc}.
In this case, the idea is to compute the Hausdorff dimension using  suitable stratifications of $M$ where the discontinuities of the dimensions $q\mapsto\dim\bD^i_q$ are somehow controlled. Namely, as suggested in \cite{gromovcc}, we consider stratifications made by submanifolds $N$ which are {\it strongly equiregular}, i.e., for which both the dimensions $\dim\bD^i_q$ and $\dim(\bD^i_q\cap T_qN)$ are constant as $q$ varies in $N$.

The first part of the paper provides an answer to question 3 when $N$ is  strongly equiregular.
The first result of the paper (Theorem~\ref{mainth}) computes
the Hausdorff dimension of a strongly equiregular submanifold  $N$ in  terms of the dimensions of $\dim(\bD^i_q\cap T_qN)$, generalizing formula \r{eqintro} which corresponds to the case $N=M$. More precisely, $\dim_H(N)=Q_N$ where
 $$
Q_N:=\sum_{i=1}^{r(p)}i(\dim(\bD^i_p\cap T_pN)-\dim(\bD^{i-1}_p\cap T_pN)).
 $$
This actually follows from a stronger property: indeed, we show that the $Q_N$-dimensional spherical Hausdorff measure in $N$  is absolutely continuous with respect to any smooth measure (i.e. any measure induced locally by a volume form) on $N$. The Radon--Nikodym derivative computed in Theorem~\ref{mainth} generalizes \cite[Lemma~32]{balu}, which corresponds to the case $N=M$.
 %As for question 2, we state Proposition~\ref{th:nuq}, which implies that $\hh^{Q_N}(B(p,\rho)\cap N)$ is positive and finite.
The main ingredient behind the proofs of such results is the fact that for a strongly equiregular submanifold $N$ the metric tangent cone to $(N,d|_N)$ exists at every $p\in N$ and can be identified to $T_pN$ via suitable systems of privileged coordinates (see Lemma~\ref{densities}).

The results for strongly equiregular submanifolds provide a first step towards the answer of questions 1 and 2 in the general case, at least for analytic sub-Riemannian manifolds. This is the topic in the second part of the paper.  Indeed, when $(M,\bD,g)$ is analytic, $M$ can be stratified as $M=\cup_{i\geq 0}M_i$ where each $M_i$ is an analytic equiregular submanifold. Then, the Hausdorff dimension of a small ball $B$ is the maximum of the Hausdorff dimensions of the intersections $B\cap M_i$. To compute the latter ones, we use that  each strata $M_i$ can further be decomposed as the disjointed union of strongly equiregular analytic submanifolds.
In Lemma~\ref{le:equireg}, using Theorem~\ref{mainth} we compute  the Hausdorff dimension of an equiregular (but possibly not strongly equiregular) analytic submanifold and we estimate the density of the corresponding Hausdorff measure.
Characterizing the finiteness of the corresponding Hausdorff measure of the intersection of a small ball with an equiregular analytic submanifold is rather involved. Yet this is the main issue in question 2, as whenever the Hausdorff measure of $B(p,\rho)\cap\{ \textrm{regular points}\}$ is infinite at a singular point $p$ then  so is $\hh^{\dim_H(M)}(B(p,\rho))$.
%Roughly speaking, if the set of singular points $\Sigma$ is strongly equiregular,   this is a consequence of  a compensation between the fact that the dimensions of the subspaces $\bD^i_q$ grow as $q$   jumps out of $\Sigma$ and the dimensions of $\bD^i_q\cap T_q\Sigma$ are constant on $\Sigma$.
To estimate $\hh^{\dim_H(M)}(B(p,\rho)\cap\{ \textrm{regular points}\})$, we assume that the singular point $p$ is ``typical'', that is, it belongs to a strongly equiregular submanifold $N$ of the singular set.
In Theorem~\ref{th:finite} we characterize the finiteness of the aforementioned measure at typical singular points through an algebraic relation involving the Hausdorff dimension $Q_{\mathrm{reg}}$ near a regular point, the Hausdorff dimension $Q_N$ of $N$, and the nonholonomic order at $p$ of the volume form on $M$ evaluated along some families of vector fields, given by Lie brackets between generators of the distribution.

The proof of Theorem~\ref{th:finite} (and of Proposition~\ref{th:nuq}) will appear in a forthcoming paper.

%\red{open questions:drop assumption A1 in Theorem~\ref{th:finite} and hence answer completely in the analytic case. Analyze whether in section 3 it is essential that the flag of $\bD$ is regular or it suffices that so is the one restricted to $N$. Answer completely in the smooth generic case}
The structure of the paper is the following. In Section~\ref{secnot} we recall shortly the definitions of Hausdorff measures and dimension and some basic notions in sub-Riemannian geometry.
Section~\ref{seceq} is devoted to the the definition and the study of strongly equiregular submanifolds and contains the proof of Theorem~\ref{mainth} and the statement of Proposition~\ref{th:nuq}. In Section~\ref{secan} we treat analytic sub-Riemannian manifolds. First, we estimate the Hausdorff dimension  $\bar Q_N$ of an analytic equiregular submanifold $N$ in Section~\ref{secand}.  Then, in Section~\ref{secanv}, we prove that the ${\bar Q}_N$-dimensional Hausdorff measure of the intersection of a small ball  $B(p,\rho)$ with $N$ is finite if $p\in N$ and we state Theorem~\ref{th:finite}. Finally, we end by applying our results to  some examples of sub-Riemannian manifolds in Section~\ref{secanex}. In particular, the examples show that when the Hausdorff dimension of a ball centered at a singular point is equal to the Hausdorff dimension of the whole manifold, the corresponding Hausdorff measure can be both finite or infinite.

\section{Basic notations}
\label{secnot}

\subsection{Hausdorff measures}

Let $(M,d)$ be
a metric space. We denote by $\diam S$ the diameter of a set $S
\subset M$, by  $B(p,\rho)$ the open ball $\{q\in M\mid
d(q,p)<\rho\}$, and by $\overline{B(p,\rho)}$ the closure of
$B(p,\rho)$.
Let $\alpha \geq 0$ be a real number. For every set $A \subset M$, the
\emph{$\alpha$-dimensional Hausdorff measure} $\hh^\alpha$ of $A$  is
defined as  $\hh^\alpha(A)
= \lim_{\eps \to 0^+} \hh^\alpha_\eps(A)$, where
$$
\hh^\alpha_\eps(A) = \inf \left\{ \sum_{i=1}^\infty  \left(\diam S_i\right)^\alpha
\, : \, A \subset \bigcup_{i=1}^\infty S_i, \ S_i \hbox{ closed set} ,
\ \diam S_i \leq \eps \right\},
$$
and the \emph{$\alpha$-dimensional spherical Hausdorff measure} is
defined as $\ss^\alpha(A)
= \lim_{\eps \to 0^+} \ss^\alpha_\eps(A)$, where
$$
\ss^\alpha_\eps(A) = \inf \left\{ \sum_{i=1}^\infty  \left(\diam
S_i\right)^\alpha \, : \, A \subset \bigcup_{i=1}^\infty S_i, \ S_i \hbox{ is
  a ball}, \ \diam
S_i \leq \eps  \right\}.
$$
For every set $A\subset M$, the non-negative number
$$
D=\sup\{\alpha\geq0\mid \hh^\alpha(A)=\infty\}=\inf\{\alpha\geq0\mid \hh^\alpha(A)=0\}
$$
is called the {\it Hausdorff dimension of $A$}. The $D$-dimensional
Hausdorff measure  $\hh^D(A)$ is called the Hausdorff volume of
$A$. Notice that this volume
may be $0$, $>0$, or $\infty$.\medskip

Given a subset $N\subset M$, we can consider the metric space $(N,d|_N)$. Denoting by $\hh^\alpha_N$ and $\ss^\alpha_N$ the Hausdorff and spherical Hausdorff measures in this space, by definition we have
\begin{eqnarray}
\hh^\alpha\llcorner_N(A)&:=&\hh^\alpha(A\cap N)=\hh_N^\alpha(A\cap N),\nonumber \\
\ss^\alpha\llcorner_N(A)&:=&\ss^\alpha(A\cap N)\leq\ss_N^\alpha(A\cap N)\label{eqes}.
\end{eqnarray}
These are a simple consequences of the fact that a  set $C$ is closed in $N$ if and only if $C=C'\cap N$, with $C'$ closed in $M$. Notice that the inequality \r{eqes} is strict in general, as coverings in the definition of $\ss^\alpha_N$ are made with sets $B$ which satisfy $B=\overline{B(p,\rho)}\cap N$ with $p\in N$, whereas coverings in the definition of $\ss^\alpha\llcorner_N$ include sets of the type  $\overline{B(p,\rho)}\cap N$ with $p\notin N$.  Moreover, by construction of Hausdorff measures, for every subset $S \subset N$,
$
\hh^\alpha(S)\leq\ss^\alpha(S)\leq 2^\alpha\hh^\alpha(S)
$
and $\hh_N^\alpha(S)\leq\ss_N^\alpha(S)\leq 2^\alpha\hh_N^\alpha(S)$.
Hence
$$
\hh^\alpha(S)\leq\ss_N^\alpha(S)\leq 2^\alpha\hh^\alpha(S),
$$
and $\ss^\alpha_N$ is absolutely continuous with respect to $\hh^\alpha\llcorner_N$.

\subsection{Sub-Riemannian manifolds} A {\it sub-Riemannian manifold
  of class $\con^k$} ($k=\infty$ or $k=\omega$ in the analytic case)
is a triplet  $(M,\bD,g)$, where $M$ is a $\con^k$-manifold,    $\bD$
is a Lie-bracket generating $\con^k$-subbundle of $TM$ of rank $m<\dim
M$ and $g$ is a Riemannian metric of class $\con^k$ on $\bD$.
Using the Riemannian metric, the length of  horizontal curves, i.e.,
absolutely continuous curves which are almost everywhere tangent to
$\bD$, is well-defined. The Lie-bracket generating assumption implies
that the  distance $d$  defined as the infimum of length of horizontal
curves  between two given points is finite and continuous
(Rashewski--Chow Theorem). We refer to  $d$ as the {\it sub-Riemannian
  distance}. The set $M$ endowed with the sub-Riemannian distance $d$
is a metric space $(M,d)$ (often called \emph{Carnot-Carath\'{e}odory
space}) which has the same topology than the
manifold $M$.

We denote by $\bD_q\subset T_qM$ the fiber of $\bD$ over $q$. The
subbundle $\bD$ can be identified with the module of sections
$$
\{X\in \VecM\mid X(q)\in\bD_q, \forall\,q\in M \}.
$$
Given $i\geq 1$, define recursively  the
  submodule $\bD^i$ of $\VecM$ by
  $$
\bD^{1}=\bD,\quad \bD^{i+1}=\bD^{i}+[\bD,\bD^{i}].
$$
Set $\bD^i_q=\{X(q)\mid X\in\bD^i\}$. Notice that the identification
between the submodule $\bD^i$ and the distribution $q\mapsto \bD^i_q$
is no more meaningful when the dimension of $\bD^i_q$ varies as a
function of $q$ (see the discussion in \cite[page 48]{bellaiche}).
The Lie-bracket generating assumption implies that for every $q\in M$
there exists an integer $r(q)$, the {\it non-holonomy degree at $q$},
such that
 \begin{equation}\label{flagd}
\{0\}\subset\bD_q^1\subset\dots\subset\bD_q^{r(q)}=T_qM.
\end{equation}
The sequence of subspaces \r{flagd} is called the {\it flag of $\bD$ at $q$}.
Set $n_i(q)=\dim\bD^i_q$ and
\begin{equation}\label{defq}
Q(q)=\sum_{i=1}^{r(q)}i(n_{i}(q)-n_{i-1}(q)),
\end{equation}
 where $n_0(q)=0$.

We say that a point $p$ is {\it regular} if, for every $i$, $n_i(q)$ is constant  as $q$ varies in a neighborhood of $p$. Otherwise, the point is said to be {\it singular}.
A subset  $A\subset M$ is called {\it equiregular} if, for every $i$, $n_i(q)$ is constant  as $q$ varies in $A$.
When the whole manifold is equiregular, the integer $Q(q)$ defined in
\r{defq} does not depend on $q$ and it is the Hausdorff dimension of
$(M,d)$ (see \cite{mitchell}).\medskip

 Given $p\in M$, let $X_1,\dots, X_m$ be a local orthonormal frame of
 $\bD$. A multiindex $I$ of length $|I|=j\geq 1$ is an element of
 $\{1,\dots,m\}^j$. With any multiindex $I=(i_1,\dots,i_j)$ is
 associated an iterated Lie bracket $X_I=[X_{i_1},[X_{i_2},
 \dots,X_{i_j}]\dots]$ (we set $X_I=X_{i_1}$ if $j=1$). The set of
 vector fields $X_I$ such that $|I|\leq j$ is a family of generators
 of the module $\bD^j$. As a consequence, if the values of
 $X_{I_1},\dots,X_{I_n}$ at $q\in M$ are linearly independent, then
 $\sum_i |I_i| \geq Q(q)$.

Let $Y$ be a vector field. We define the  {\it length of $\,Y$} by
$$
\ell(Y)=\min\{i\in \N\mid Y\in\bD^i\}.
$$
In particular, $\ell (X_I) \leq |I|$. Note that, in general, if a vector field $Y$ satisfies $Y(q)\in\bD^i_q$ for every $q\in M$,  $Y$ need not be in the submodule $\bD^i$.
By an {\it adapted basis} to the flag \r{flagd} at $q$, we mean $n$ vector fields $Y_1,\dots, Y_n$  such that their values at $q$ satisfy
$$
\bD^i_q=\span\{Y_j(q)\mid  \ell(Y_j)\leq i\},\quad \forall\, i=1,\dots, r(q).
$$
In particular, $\sum_{i=1}^{n}\ell(Y_i) = Q(q)$. As a consequence, a family of brackets $X_{I_1},\dots,X_{I_n}$ such that $X_{I_1}(q),\dots,X_{I_n}(q)$ are linearly independent is an adapted basis to the flag \r{flagd} at $q$ if and only if $\sum_i |I_i| = Q(q)$.

\section{Hausdorff dimensions and volumes of
strongly equiregular   submanifolds}
\label{seceq}

In this section, we answer question 3 when $N$ is a particular kind of submanifold, namely a strongly equiregular one.  These
results include the case where $M$ itself is equiregular.

\subsection{Strongly equiregular submanifolds}

Let $N\subset M$ be a smooth connected submanifold of dimension
$b$. The {\it flag  at $q\in N$ of $\bD$ restricted to $N$} is the sequence of subspaces
\begin{equation} \label{flagq}
\{0\}\subset(\bD_q^1\cap T_qN)\subset\dots\subset(\bD_q^{r(q)}\cap T_qN)=T_qN.
\end{equation}
Set
$$%\begin{eqnarray*}
n_i^N(q)= \dim (\bD_q^i\cap T_qN)\qquad \hbox{and} \qquad
Q_N(q)=\sum_{i=1}^{r(q)}i(n^N_{i}(q)-n^N_{i-1}(q)),%\label{defqn}
$$%\end{eqnarray*}
with $n^N_{0}(q)=0$.
\begin{definition}\label{sequi}
 We say that
$N$ is {\it strongly equiregular} if
\bi
\iii[(i)] ~~$N$ is equiregular, that is, for every $i$,  the dimension
 $n_i(q)$ is constant as $q$ varies in $N$.
\iii[(ii)]  ~~for every $i$,  the dimension
 $n_i^N(q)$ is constant as $q$ varies in $N$.
 \ei
 In this case, we denote by $Q_N$  the constant value of $Q_N(q)$,
 $q\in N$.
 \end{definition}

By an {\it adapted basis} to the flag \r{flagq} at $q\in N$, we mean $b$ vector fields $Z_1,\dots, Z_b$ such taht
$$
\bD^i_q\cap T_qN=\span\{Z_j(q)\mid   \ell(Z_j)\leq i \},\quad \forall\, i=1,\dots, r(q).
$$
In particular, when $Z_1,\dots, Z_b$ is adapted to the flag \r{flagq},
we have $T_qN=\span\{Z_1(q),\dots, Z_b(q)\}$ and
$\qn=\sum_{i=1}^{b}\ell(Z_i)$. \medskip

% \red{ meaning of this definition: strong equiregularity implies
% that, on one hand, we can find privileged coordinates which rectify
% the whole manifold $N$, on the other hand, these coordinates can be
% chosen smoothly as $q$ varies in $N$.}

\medskip

 Recall that the metric tangent cone\footnote{in Gromov's sense, see \cite{gromov}} to $(M,d)$ at any point $p$ exists and it is isometric to $(T_pM,\widehat d_p)$, where  $\widehat d_p$ denotes the sub-Riemannian distance associated with a nilpotent approximation at $p$ (see \cite{bellaiche}).
The following lemma shows the relevance of strongly equiregular submanifolds as particular subsets of $M$ for which a metric tangent cone  exists. Such metric space is isometrically embedded in a metric tangent cone to the whole $M$ at the point.
%Moreover,  the   smooth volume of balls centered on the submanifold can be estimated in terms of the induced volume in the metric tangent cone.

\bl\label{densities}
Let $N\subset M$ be a $b$-dimensional submanifold of $M$. Assume $N$
is strongly equiregular.
%\footnote{For points (i), (ii) and the first part of point (iii)  it
%suffices to assume that the dimensions
%$n_i^N(q)$ are constant   $q$ varies in $N$. For the second part of
%point (iii)  we still need to check that convergence in \r{measure}
%is uniform.}
Then, for every $p\in N$:
\bi
\item[$\mathrm{(i)}$]  ~there exists a metric tangent cone to
  $(N,d|_N)$ at $p$ and it is isometric to $(T_pN,\widehat
  d_p|_{T_pN})$;
\item[$\mathrm{(ii)}$]  ~the graded vector space
 $$
 \mathfrak{gr}_p^N(\bD):=\oplus_{i=1}^{r(p)}(\bD^{i}_p\cap T_pN)/(\bD^{i-1}_p\cap T_pN)
 $$
 is a nilpotent Lie algebra whose associated Lie group $\mathrm{Gr}_p^N(\bD)$ is diffeomorphic to $T_pN$;
 \iii[$\mathrm{(iii)}$]  ~~every $b$-form $\omega\in\bigwedge^b N$ on
 $N$ induces canonically a  left-invariant $b$-form $\wp$ on
 $\mathrm{Gr}_p^N(\bD)$. Moreover,
\begin{equation}\label{measure}
\int_{N\cap B(p,\eps)}\omega =\eps^{Q_N}\int_{T_pN\cap \widehat B_p}\wp+o(\eps^{Q_N}),
\end{equation}
where $o(\eps^{Q_N})$ is uniform as $p$ varies in $N$ and $\widehat
B_p$ is the ball centered at $0$ of radius $1$ in the nilpotent
approximation at $p$ of the sub-Riemannian manifold.
\ei
\el
\brem
When $N$ is an open submanifold of $M$, assuming $N$ strongly
equiregular is equivalent to saying that $N$ contains only regular
points. In that case,   Lemma~\ref{densities} is well-known (point (i)
follows by the fact that  the nilpotent approximation is a metric
tangent cone, point (ii) says that the tangent cone shares a group
structure - which in this case satisfies the additional property
$\mathfrak{gr}_p(\bD)=\mathrm{span}_p\{\bD^1\}$ - and (iii) has been
remarked in \cite{balu} using the canonical isomorphism between
$\bigwedge^n(\mathfrak{gr}_p(\bD)^*)$ and $\bigwedge^n(T^*_pM)$.
\erem

{\it Proof.}
Note first that since the result is of local nature, it is sufficient
that we prove it on a small neighbourhood $B(p_0,\rho)\cap N$ of a point
$p_0 \in N$. For every $p$ in a such a neighbourhood, there exists a
coordinate system
$\varphi_p:U_p\to \R^n$ on a neighborhood $U_p \subset M$ of $p$, such that
 $\varphi_p$ are privileged coordinates at $p$,  $p \mapsto \varphi_p$ is continuous, and $N$ is rectified in coordinates $\varphi_p$, that is $\varphi_p(N \cap U_p) \subset \{x\in\R^n\mid x_{b+1}=\dots=x_n=0\}$.
The construction is as follows.

Given $\rho>0$ small enough, we can find   $b$ vector fields
$Y_1,\dots, Y_b$ defined on $B(p_0,\rho)$ which form a
basis adapted to the flag \r{flagq}  restricted to $N$ at every $p\in
B(p_0,\rho)\cap N$. Moreover, up to reducing $\rho$, we can find
$Y_{b+1},\dots, Y_{n}$ such that $Y_1,\dots, Y_n$ is adapted to the
flag \r{flagd} of the distribution at every point $p \in B(p_0,\rho)\cap
N$.
Using these bases, we define
for   $p \in N\cap B(p_0,\rho)$, a local diffeomorphism
$\Phi_p:\R^n\to M$ by
\begin{equation}\label{coorn}
\Phi_p(x)=\exp\left(\sum_{i=b+1}^nx_iY_i\right)\circ\exp
\left(\sum_{i=1}^bx_i Y_i\right)(p).
\end{equation}
 The inverse $\varphi_p=\Phi_p^{-1}$ of $\Phi_p$ provides a system of
 coordinates centered at
 $p$ which are privileged (see \cite{hermes}). Moreover, thanks to
 property (i) in Definition~\ref{sequi},  the map from $B(p_0,\rho) \cap
 N$ to $M$ which associates with $p$ the point $\Phi_p(x)$ is smooth
 for every $x\in\R^n$. Finally, in coordinates $\varphi_p$, the submanifold
 $N \cap U$ coincides with the set
 $$
 \left\{\exp\left(\sum_{i=1}^bx_i Y_i\right)(p)\mid
   (x_1,\dots,x_b)\in \Omega \right\} \subset \left\{ \Phi_p(x)\mid
   x_{b+1}=\dots=x_n=0 \right\},
 $$
where $\Omega$ is an open subset of $\R^b$.

Using $\varphi_p$ we identify $M$ with $T_pM\simeq\R^n$. Since
$Y_1(p),\dots, Y_b(p)$ span $T_pN$, $\varphi_p$ maps $N$ in $T_pN$,
where $T_pN$ is identified
with  $\R^b\times\{0\}\subset \R^n\simeq T_pM$.
Therefore, whenever $q_1,q_2\in U\cap N$ we have
$$
\widehat d_p(q_1,q_2)=\widehat d_p|_{T_pN}(q_1,q_2),
$$
and obviously $d(q_1,q_2)= d|_{N}(q_1,q_2)$. Hence estimate (70) in
\cite[Theorem~7.32]{bellaiche} holds when we restrict $d$ to $N$ and
$\widehat d$ to $T_qN$. This allows to conclude that a metric tangent
cone to $(N,d|_N)$ at $p$ exists and it is isometric to
$(T_pN,\widehat d_p|_{T_pN})$, where the inclusion of $T_pN$ into
$T_pM$ is to be intended via $\varphi_p$.

The algebraic structure of $\mathfrak{gr}_p^N(\bD)$ and the  fact that
$Gr_p^N(\bD)$ is diffeomorphic to $\R^b$ are straightforward. As a
consequence, there also exists a canonical isomorphism between
$\bigwedge^b(\mathfrak{gr}^N_p(\bD)^*)$ and $\bigwedge^b(T^*_pN)$. Let
$\tilde\omega_p$ be the image of $\omega_p$ under such isomorphism
(see the construction in \cite[Section~10.5]{montgomery}). Then $\wp$
is defined as the left-invariant $b$-form on $T_pN$ which coincides
with $\tilde\omega_p$ at the origin.

Finally, as a consequence of point (i),  by definition of metric tangent cone
$\varphi_p(B(p,\eps)\cap N)$ converges to $\widehat B(0,\eps)\cap
T_pN$ in the Gromov--Hausdorff sense as $\eps$ goes to $0$. By
homogeneity of $\widehat d_p$ we have
$\widehat B(0,\eps)\cap T_pN=\eps^{Q_N}(\widehat B_p\cap T_pN)$ and we
get \r{measure}. Since  $p\mapsto\varphi_p$ and $p\mapsto \widehat
B_p$ are  continuous \cite[Section~4.1]{balu},   the
remainder $o(\eps^{Q_N})$ in \r{measure} is uniform with respect to
$p$.
\hfill$\square$

\medskip

For the sake of completeness, let us give an explicit formula for
$\wp$. Recall that the
construction of the coordinates $\varphi_p$ involves an adapted basis
$Y_1,\dots, Y_b$  to the flag \r{flagq}  restricted to $N$ at every
$p\in B(p_0,\rho)\cap N$. In particular the vector fields $Y_1,\dots,
Y_b$ restricted to $N$ form a local frame for the tangent bundle to
$N$ and
$$
\omega =\omega(Y_1,\dots,Y_b) d( Y_1|_{N})\wedge\dots\wedge d( Y_b|_{N}).
$$
  Let $X_1,\dots, X_m$ be a local orthonormal frame for the sub-Riemannian structure in a neighborhood of $p$, and  $X_{I_1},\dots, X_{I_n}$ be an adapted basis to the flag \r{flagd} at $p$, where  $X_{I_j}$ is the Lie bracket corresponding to the multi-index $I_j$.
Since $X_{I_1},\dots, X_{I_n}$  is a local frame for the tangent bundle to $M$, for every $i=1,\dots, b$ we can write $Y_i$ in this basis  as
$$
Y_i=\sum_{|I|\leq \ell(Y_i)}Y_i^I X_{I},
$$
where $Y_i^I$ are smooth function (the fact that only multiindices with length smaller than $\ell(Y_i)$ appear in this sum is due to the definition of length of a vector field).
Denote by $\widehat X_1^p,\dots, \widehat X_m^p$ the nilpotent approximation of $X_1,\dots, X_m$ at $p$ obtained in coordinates $\varphi_p$, and by $\widehat X^p_{I_j}$ the Lie bracket between the $\widehat X_1^p,\dots, \widehat X_m^p$ corresponding to the multiindex $I_j$.
For every $i=1,\dots, b$ we define the vector field
$$
\widehat Y^p_i=\sum_{|I|= \ell(Y_i)}Y_i^I(p) \widehat X_{I}.
$$
This enables us to compute $\wp$ as
\begin{equation}\label{eqwp}
\wp =\omega_p(Y_1(p),\dots,Y_b(p))d(\widehat Y^p_1|_{T_pN})\wedge\dots\wedge d(\widehat Y^p_b|_{T_pN}).
\end{equation}
The fact that the right-hand side of \r{eqwp}  does not depend on the
$X_I$ nor on the $Y_i$ is a consequence of the intrinsic definition of
$\wp$.%\medskip

\subsection{Hausdorff volume}
Assume now that $N$ is an orientable submanifold. By a {\it smooth
  volume} on $N$ we
mean a measure $\mu$ associated with a never vanishing  smooth form
$\omega\in\bigwedge^bN$,
 i.e., for every Borel set $A\subset N$, $\mu(A)=\int_A\omega$. We
 will denote by $\mup$ the smooth volume on $T_pN$ associated  with
 $\wp$.
% As the results in this section are local, for every $p\in N$, when
% $\rho$ is small enough, $N\cap B(p,\rho)$ is orientable.

We are now in a position to prove the main result.

\bt\label{mainth}
Let $N\subset M$ be a smooth orientable submanifold. Assume $N$ is strongly equiregular. Then, for every smooth volume $\mu$ on $N$,
\begin{equation}\label{abscont}
\lim_{\eps\to 0}\frac{\ss^{\qn}_N(B(q,\eps))}{\mu(N \cap
  B(q,\eps))}=\frac{\diam_{\widehat d_q}(T_qN\cap \widehat
  B_q)^{Q_N}}{\muq(T_qN\cap \widehat B_q)},~~ \forall\, q\in N,
\end{equation}
where $\diam_{\widehat d_q}$ denotes the diameter with respect to the distance $\widehat d_q$.
In particular, $\ss^{\qn}_N$ is absolutely continuous with respect to $\mu$ with Radon--Nikodym derivative   equal to the right hand side of \r{abscont}.
As a consequence,
\begin{equation}\label{dimh}
\dim_{H}N=\qn,
\end{equation}
and, for a small ball $B(p,\rho)$ centered at a point $p\in N$, the
Hausdorff volume $\hh^{\qn}(N \cap B(p,\rho))$ is finite.
\et

\brem\label{regfin}
When $N$ is an open submanifold of $M$, e.g., $N=\{p\in M\mid p \mbox{ is regular}\}$, the computation of Hausdorff dimension is well-known, see  \cite{mitchell}. In particular, when $p$ is a regular point the top-dimensional Hausdorff measure $\hh^{Q}(B(p,r))$ is positive and finite. When $N=M$, equation \r{abscont} gives a new proof to \cite[Theorem~1]{balu}. This is interesting since the latter was obtained as a consequence of \cite[Lemma~32]{balu}, whose proof is incorrect.
%Since the  proof of such Lemma  is incorrect, we find it important to include here the proof of Theorem~\ref{mainth}.}
\erem
To prove Theorem~\ref{mainth} a fundamental step is the following lemma.

\bl
\label{le:ss_vs_mu}
Let $N$ and $\mu$ be as in Theorem~\ref{mainth}. Let $p\in N$. Assume   there exists positive constants $\eps_0$ and $\mu_+ > \mu_-$ such that, for every $\eps<\eps_0$ and every point $q \in B(p,\eps_0)\cap N$, there holds
\begin{equation}
\label{eq:mu+-}
    \mu_- \diam(B(q,\eps)\cap N)^{Q_N} \leq \mu (B(q,\eps)\cap N) \leq \mu_+ \diam(B(q,\eps)\cap N)^{Q_N}.
\end{equation}
Then, for every $\eps<\eps_0$,
$$
\frac{\mu (B(p,\eps)\cap N)}{\mu_+}  \leq \ss_N^{Q_N}(B(p,\eps)) \leq  \frac{\mu (B(p,\eps)\cap N)}{\mu_-}.
$$
\el

\noindent{\it Proof.}
Let $\bigcup_i B(q_i,r_i)$ be a covering of $B(p,\eps)\cap N$ with
balls of radius smaller than $\delta <\eps_0$. If $\delta$ is small enough,
every $q_i$ belongs to $B(p,\eps_0)\cap N$ and, using~\eqref{eq:mu+-}, there holds
$$
 \mu(B(p,\eps)\cap N) \leq \sum_i \mu (B(q_i,r_i)\cap N) \leq \mu_+
\sum_i \diam (B(q_i,r_i)\cap N)^{Q_N}.
$$
Hence, we have $\ss_N^{Q_N} (B(p,\eps)) \geq \frac{\mu(B(p,\eps)\cap N)}{\mu_+}$.

For the other inequality, let $\eta>0$, $0<\delta <\eps_0 $ and let $\bigcup_i B(q_i,r_i)$ be a covering of $B(p,\eps)\cap N$ such that $q_i\in B(p,\eps)\cap N$ $r_i <\delta$ and $\sum_i \mu (B(q_i,r_i)\cap N) \leq \mu(B(p,\eps)) + \eta$. Such a covering exists due to the Vitali covering lemma. Using as above~\eqref{eq:mu+-}, we obtain
$$
\mu(B(p,\eps)\cap N) + \eta \geq \sum_i \mu (B(q_i,r_i)\cap N) \geq \mu_-
\sum_i \diam(B(q_i,r_i)\cap N)^{Q_N}.
$$
We then have $\ss_{N,\delta}^\qn (B(p,\eps))  \leq  \frac{\mu(B(p,\eps)\cap N)}{\mu_-}  + \frac{\eta}{\mu_-} $. Letting  $\eta$ and $\delta$ tend to $0$, we get the conclusion.
\hfill$\square$

\bigskip

\noindent{\it Proof of Theorem~\ref{mainth}.} Fix $q\in N$.  By point $(ii)$ in Lemma~\ref{densities} $(T_qN,\widehat d_q|_{T_qN})$ is a metric tangent cone to $(N,d|_N)$ at $q$, whence, from the definition of Gromov--Hausdorff convergence we get
\begin{equation}\label{limdiam}
\lim_{\eps\to 0}\frac{\diam(N\cap B(q,\eps))}{\eps}=\diam\,_{\widehat d_q}(T_qN\cap \widehat B_q).
\end{equation}
By \r{measure} in Lemma~\ref{densities}, for every $q\in N$ there holds
\begin{equation}\label{limmu}
\mu(N\cap B(q,\eps)) = \eps^{Q_N} \muq (T_qN\cap \widehat B_q) + o(\eps^{Q_N}).
\end{equation}
Since $N$ is strongly equiregular, the limits in \r{limdiam} and \r{limmu} hold uniformly as $q$ varies in $N$.

 Moreover, adapting the argument in \cite[Section~4.1]{balu}, we deduce that the map $q\mapsto  \muq (\widehat B_q\cap T_qN)$ is continuous on $N$.
 As a consequence, for any $\eta>0$ there exists $\eps_1>0$ such that for every $q\in B(p,\eps_1)$ and every $\eps<\eps_1$ we have
$$
\mu_-\leq\frac{\mu(N\cap B(q,\eps))}{\diam(N\cap B(q,\eps))^\qn}\leq\mu_+
$$
with
$$
\mu_\pm = \frac{\widehat \mu_q(T_qN\cap \widehat B_q)}{\diam_{\widehat d_q}(T_qN\cap \widehat B_q)^{Q_N}}\pm \eta.
$$
Therefore, applying Lemma~\ref{le:ss_vs_mu} and letting $\eta$ tend to $0$ we deduce \r{abscont}.

To show \r{dimh}, notice that  the right-hand side of \r{abscont} is
continuous and positive as a function of $q$. Hence,
 for $\ss^{\qn}_N$-almost every $q\in N$ there exists $\rho>0$ small
 enough such that
\begin{equation}\label{esq}
0<\ss^{\qn}(N\cap B(p,\rho))<\infty.
\end{equation}
This is equivalent to \r{dimh}.
\hfill$\square$

\bigskip

We end this section by stating a result which gives a weak
  equivalent of the function $\muq(T_qN\cap \widehat B_q)$ appearing
  in Theorem~\ref{mainth}. This will be useful in the following to
  determine whether the Hausdorff volume of a small ball is finite or
  not. This result stems from the uniform Ball-Box Theorem,
  \cite{jean1} and \cite[Th.\ 4.7]{cimpa}.

\begin{proposition}
\label{th:nuq}
Let $M$ be orientable and $\varpi$ be a volume form on $M$.
Let $N$ be an orientable submanifold of $M$ of dimension $b$, and let
$\omega$ be a volume form on $N$, with associated smooth volume
$\mu$. Assume $N$ is strongly equiregular and set $Q[N]$ equal to the
constant value of $Q(q)$, for $q\in N$. Then there exists a constant
$C>0$ such that, for every $q \in N$,
$$
\frac{1}{C} \nu_q \leq \muq(T_qN\cap \widehat B_q) \leq C \nu_q \quad
(\hbox{i.e. } \muq(T_qN\cap \widehat B_q) \asymp \nu_q \hbox{
  uniformly w.r.t.\ $q$}),
$$
where $\nu_q = \max\{ \left(\omega\wedge dX_{I_{b+1}} \wedge \cdots
  \wedge dX_{I_{n}}\right)_q (X_{I_{1}}(q),\dots,X_{I_{n}}(q))\}$, the
maximum being taken among all $n$-tuples $(X_{I_{1}},\dots,X_{I_{n}})$
in $\arg \max \{ \varpi_q(X_{I'_{1}}(q),\dots,X_{I'_{n}}(q)) \mid \
\sum_i |I'_i| = Q[N] \}$.

In particular, if $N$ is an open equiregular subset of $M$, i.e.,
$b=n$, and if $\mu$ is the smooth measure on $M$ associated with
$\varpi$, we have
$$
\muq(\widehat B_q) \asymp \max \{
\varpi_q(X_{I'_{1}},\dots,X_{I'_{n}}) \mid \ \sum_i |I'_i| = Q [M] \},
\quad \hbox{ uniformly w.r.t.\ $q \in M$}.
$$
\end{proposition}

This proposition, together with Theorem~\ref{mainth}, allows to give an estimate of the Hausdorff volume of a subset of $N$. If $S \subset N$, then
\begin{equation}\label{eq:asymhh}
\frac{1}{C'} \int_S \frac{1}{\nu_q} d\mu \leq \hh^{Q_N} (S) \leq C' \int_S \frac{1}{\nu_q} d\mu,
\end{equation}
where the constant $C'>0$ does not depend on $S$.

\smallskip

\section{Hausdorff dimensions and volumes of analytic
sub-Riemannian manifolds}\label{secan}

%%%%%%%%%%%%%%%%%%%%%%%%
% \red{Questions:
% \begin{itemize}
%   \item Ajouter def de Whitney stratification?? Yes, also because at
%   a certain point you say that by definition a stratification is
%   locally finite, which in principle is not granted.
% \end{itemize}}

%%%%%%%%%%%%%%%%%%%%%%%%

Let $(M,\bD,g)$ be an analytic ($C^\omega$) sub-Riemannian manifold. The set $\Sigma$ of singular points is an analytic subset of $M$ which admits a Whitney stratification $\Sigma=\bigcup_{i\geq 1} M_i$ by analytic and equiregular submanifolds $M_i$ (see for instance~\cite{Goresky1988}).
Denoting $M_0=M \setminus \Sigma$ the set of regular points, we obtain a Whitney stratification $M=\bigcup_{i\geq 0} M_i$ of $M$ by analytic and equiregular submanifolds. Note that $M_0$ is an open and dense subset of $M$, but  it may be disconnected. As a consequence, the Hausdorff dimension of $M$ is
$$
\dim_H( M)=\max_{i\geq 0} \dim_H(M_i),
$$
and the $\alpha$-dimensional Hausdorff measure of a ball $B(p,\rho)$, $p\in M$ and $\rho>0$, is
$$
\hh^{\alpha}(B(p,\rho))=\sum_i\hh^{\alpha}(B(p,\rho)\cap M_i).
$$

\subsection{Hausdorff dimension}\label{secand}
The first problem is then to determine the Hausdorff dimension of an equiregular - possibly not strongly equiregular - submanifold.

\bl
\label{le:equireg}
Let $N$ be an analytic and equiregular submanifold of $M$. Set $\qon := \max_{q \in N} \qn(q)$. Then
$$
\dim_H( N)= \qon,
$$
and $\qn(q) = \qon$ on an open and dense subset of $N$.

If moreover $N$ is orientable, then for every smooth measure $\mu$ on $N$, $\ss^{\qon}_N$ is absolutely continuous with respect to $\mu$ with Radon--Nikodym derivative
\begin{equation}
\label{eq:rdnk_equireg}
 \frac{d\ss^{\qon}_N}{d\mu} (q) = \frac{(\diam_{\widehat d_q}(T_qN\cap \widehat B_q))^{\qon}}{\widehat \mu_q(T_qN\cap \widehat B_q)}, \quad \hbox{for $\mu$-a.e.\ } q\in N.
\end{equation}
\el

\begin{proof}
Since $N$ is analytic and equiregular, it admits a stratification $N= \bigcup_{i} N_i$ by strongly equiregular submanifolds $N_i$ of $N$. By Theorem~\ref{mainth}, $\dim_H( N_i)=Q_{N_i}$ and thus $\dim_H( N)=\max_{i} Q_{N_i}$. In particular, $\dim_H( N) \leq \max_{q \in N} \qn (q)$.

Now, recall that $\qn(q) = \sum_{i=1}^{r_N}i(n_i^N(q)-n_{i-1}^N(q))$, where $r_N:=r(q)$ is constant since $N$ is equiregular, and $n_{r_N}^N(q) = \dim N$. This may be rewritten as
\begin{equation}
\label{eq:qn}
\qn(q) = \sum_{i=0}^{r_N-1} \mathrm{codim} (\bD_q^i\cap T_qN),
\end{equation}
where $\mathrm{codim} (\bD_q^i\cap T_qN) = n_{r_N}^N(q)-n_{i}^N(q)$ is
the codimension of $\bD_q^i\cap T_qN$ in $T_qN$. The submanifold $N$
being equiregular, $\qn(q)$ is a lower semi-continuous function on $N$
with integer values. Hence $\qn(q)$ takes its maximal value $\qon$ on
the strata $N_i$ which are open in $N$, and smaller values on non open
strata. Since $Q_{N_i}(q) = Q_N (q)$ when $N_i$ is an open subset of
$N$ and $Q_{N_i}(q) < Q_N(q)$ when $N_i$ is a non open subset of $N$,
the first part of the lemma follows.

As for the second part, notice that every non open stratum $N_i$ is of $\mu$-measure zero, since $N_i$ is a subset of $N$ of positive codimension, and of $\ss^{\qon}_N$-measure zero, since $\dim_H( N_i)=Q_{N_i}< \qon$. A first consequence is that $N$ is strongly equiregular near $\mu$-a.e.\ point $q$. Therefore the measure $\widehat \mu_q$ on $T_qN$ is defined $\mu$-a.e.\ -- and so is the right-hand side of~\eqref{eq:rdnk_equireg}. Applying then Theorem~\ref{mainth} to every open stratum $N_i$, we get the conclusion.
\hfill$\square$
\end{proof}

\begin{corollary}
\label{cor:dimHM}
$\quad  \displaystyle
\dim_H( M)=\max \{ Q_{M_i}(q) \ : \ i\geq 0, \ q \in M_i \} = \max \{ \overline{Q}_{M_i} \ : \ i\geq 0 \}
$.
\end{corollary}
\medskip

\subsection{Finiteness of the Hausdorff volume of balls}\label{secanv}

Let $p\in M$ and $\rho>0$ ($\rho$ is assumed to be arbitrarily small). The aim of this section is to determine under which conditions the small ball $B(p,\rho)$ has a finite Hausdorff volume $\hh^{\dim_H(B(p,\rho))}(B(p,\rho))$.
We make first two preliminary remarks.
\begin{itemize}
  \item If $p$ is a regular point, then there exists a neighbourhood of $p$ in $M$ which is strongly equiregular, and Theorem~\ref{mainth} implies that $\hh^{\dim_H(B(p,\rho))}(B(p,\rho))$ is finite. We then assume in the following that $p$ is a singular point.
  \item The results of this section are local. Up to reducing to a neighbourhood of $p$, we can assume that $M$ is an oriented manifold with volume form $\varpi$.
\end{itemize}

Recall that, by definition, the stratification $M=\bigcup_{i\geq 0} M_i$ is locally finite. That is, there exists a finite set $\mathcal{I}$ of indices such that $p \in \overline{M_i}$ if and only if $i \in \mathcal{I}$, where $\overline{M_i}$ denotes the closure of the stratum $M_i$. Therefore,
for $\rho$ small enough, the ball $B(p,\rho)$ admits a finite stratification $B(p,\rho)= \bigcup_{i \in \mathcal{I}} (B(p,\rho) \cap M_i)$. Applying Corollary~\ref{cor:dimHM}, the Hausdorff dimension $D_p$ of $B(p,\rho)$ is
$$
D_p = \max \{ Q_{M_i}(q) \ : \ i \in \mathcal{I}, \ q \in M_i \}.
$$
Let $\mathcal{J} \subset \mathcal{I}$ be the subset of indices $i$ such that $\dim_H(M_i)=D_p$. We have
$$
\hh^{D_p}(B(p,\rho))=\sum_{ i \in \mathcal{J}}\hh^{D_p}(B(p,\rho)\cap M_i).
$$

\bp
\label{le:pinN}
Let $N$ be an analytic and equiregular submanifold of $M$, $\dim_H( N)= \qon$. If $p\in N$ and if $\rho >0$ is small enough, then the Hausdorff volume $\hh^{\qon}(B(p,\rho)\cap N)$ is finite.
\ep

\begin{proof}%[of Proposition~\ref{le:pinN}]
Up to replacing $N$ with a small neighbourhood of $p$ in $N$, we assume that $N$ is orientable. We then choose a smooth measure $\mu$ on $N$ and we have, for $\rho$ small enough, $\mu (B(p,\rho)\cap N) < +\infty$. From Lemma~\ref{le:equireg},
$$
 \ss^{\qon}_N (B(p,\rho)\cap N) = \int_{B(p,\rho)\cap N}\frac{(\diam_{\widehat d_q}(T_qN\cap \widehat B_q))^{\qon}}{\widehat \mu_q(T_qN\cap \widehat B_q)} d\mu.
$$
The submanifold $N$ is strongly equiregular near $\mu$-a.e.\ $q\in N$. We can then apply Proposition~\ref{th:nuq} near $\mu$-a.e.\ $q\in N$ and we get
$$
 \ss^{\qon}_N (B(p,\rho)\cap N) \leq C \int_{B(p,\rho)\cap N}\frac{(\diam_{\widehat d_q}(T_qN\cap \widehat B_q))^{\qon}}{\nu_q} d\mu.
$$
The function $q \mapsto \nu_q$ is positive and continuous on $N$, so the integrand function in the previous formula is finite and continuous on $N$, and we have $\ss^{\qon}_N (B(p,\rho)\cap N) \leq \mathrm{Cst}\ \mu (B(p,\rho)\cap N) < +\infty$. Since $\hh^{\qon}$ is absolutely continuous with respect to $\ss^{\qon}_N$, the conclusion follows.
\hfill$\square$
\end{proof}

As a consequence, the Hausdorff volume $\hh^{D_p}(B(p,\rho))$ is finite if and only if $\hh^{D_p}(B(p,\rho)\cap M_i)$ is finite for every stratum $M_i$ such that $\dim_H(M_i)=D_p$ and $p \in \partial M_i$. To go further, we will assume that $p$ is a \emph{typical} singular point, that is, that $p$ satisfies the following assumptions for $\rho$ small enough:
\begin{description}
  \item[(A1)] $p$ belongs to a strongly equiregular submanifold $N$ of $M$, $N \subset \Sigma$, and $B(p,\rho) \cap \Sigma \subset N$;
  \item[(A2)] for every $q \in N \cap B(p,\rho)$, there exists a family $X_{I_{1}},\dots,X_{I_{n}}$ such that $\sum_i |I_i| = Q_{\mathrm{reg}}$ and $\mathrm{ord}_q \varpi (X_{I_{1}},\dots,X_{I_{n}}) = \sigma$, where $Q_{\mathrm{reg}}$ is the constant value of $Q(q)$ for $q \in M\setminus\Sigma$, and
  $$
  \sigma = \max \{ s \in \N \ : \  q \in N \cap B(p,\rho)  \ \hbox{and} \ \sum_i |I_i| = Q_{\mathrm{reg}} \ \hbox{imply} \ \mathrm{ord}_q \varpi (X_{I_{1}},\dots,X_{I_{n}}) \geq s \}.
  $$
\end{description}

Let us recall the definition of $\mathrm{ord}_q$ (see \cite{bellaiche} for details).  Given $f\in\con^k(M)$, we say that $f$ has  {\it non-holonomic order at $p$ greater than or equal to $s$}, and we write $\ord_p f \geq s$ if for every $j\leq s-1$
$$
(X_{i_1}\dots X_{i_j}f)(p)=0\quad \forall~ (i_1,\dots,i_j)\in\{1,\dots,m\}^j ,
$$
where  $X_i f$ denotes the Lie derivative of $f$ along $X_i$. Equivalently, $f(q) = O(d(p,q)^s)$. If moreover we do not have $\ord_p f \geq s+1$, then we say that $f$ has {\it  non-holonomic order at $p$ equal to $s$}, and we write $\ord_pf=s$.

\bt
\label{th:finite}
Assume $p$ satisfies (A1) and (A2). Let
%$Q_{\mathrm{reg}}$ be the
%constant value of $Q(q)$, $q \in M\setminus\Sigma$,
 $Q_N$ be the
constant value of $Q_N(q)$ for $q\in N$, and $r_{\not N}$ be the maximal
integer $i$ such that $n_i(p)-n_{i-1}(p) >
n^N_i(p)-n^N_{i-1}(p)$. Then
$$
\hh^{Q_{\mathrm{reg}}}(B(p,\rho) \setminus\Sigma)<\infty \quad \Leftrightarrow \quad \sigma \leq Q(p) - Q_N - r_{\not N}.
$$
As a consequence,
\begin{itemize}
  \item if $Q_{\mathrm{reg}}<Q_N$, then $D_p=Q_N$ and $\hh^{D_p}(B(p,\rho))$ is finite;
  \item if $Q_{\mathrm{reg}} \geq Q_N$, then $D_p=Q_{\mathrm{reg}}$ and $\hh^{D_p}(B(p,\rho))$ is finite if and only if $\sigma \leq Q(p) - Q_N - r_{\not N}$.
\end{itemize}
\et

The proof of this theorem is postponed to a forthcoming paper. It relies on the use of Proposition~\ref{th:nuq}.

\begin{remark}
Assumption (A2) is actually not necessary for the computations. If $p$ satisfies only (A1), we introduce two integers $\sigma_- \leq \sigma_+$:
$$
\begin{array}{lll}
  \sigma_+ & = & \min \{ s \in \N \ : \ \forall q \in N \cap B(p,\rho), \ \exists X_{I_{1}},\dots,X_{I_{n}} \ \hbox{s.t.} \ \sum_i |I_i| = Q_{\mathrm{reg}} \ \hbox{and} \ \mathrm{ord}_q \varpi (X_{I_{1}},\dots,X_{I_{n}}) \leq s \},\\
   \sigma_- & = & \max \{ s \in \N \ : \ \exists \hbox{ an open subset $\Omega$ of }  N \cap B(p,\rho) \  \hbox{s.t.} \ q \in \Omega  \ \hbox{and} \ \sum_i |I_i| = Q_{\mathrm{reg}} \ \hbox{imply} \ \mathrm{ord}_q \varpi (X_{I_{1}},\dots,X_{I_{n}}) \geq s \}.
\end{array}
$$
Assumption (A2) is equivalent to $\sigma_- = \sigma_+ = \sigma$. The generalization of the criterion of Theorem~\ref{th:finite} to the case where $p$ satisfies only (A1) is then:
\begin{itemize}
  \item if $\sigma_+ \leq Q(p) - Q_N - r_{\not N}$, then $\hh^{Q_{\mathrm{reg}}}(B(p,\rho) \setminus\Sigma)<\infty$;
  \item if $\sigma_- > Q(p) - Q_N - r_{\not N}$, then $\hh^{Q_{\mathrm{reg}}}(B(p,\rho) \setminus\Sigma) = \infty$.
\end{itemize}
\end{remark}

Notice that the order $\sigma$  (and $\sigma_-$ if $p$ does not satisfies (A2)) always satisfies $\sigma \geq Q(p)-Q_{\mathrm{reg}}$. We thus obtain a simpler criterion for the non finiteness of the Hausdorff volume of a ball.

\begin{corollary}
\label{cor:infinite}
Assume $p$ satisfies (A1). If $0\leq Q_{\mathrm{reg}}-Q_N < r_{\not N}$, then $\hh^{D_p}(B(p,\rho)) = \infty$.
\end{corollary}

\subsection{Examples}\label{secanex}

\begin{example}[the Martinet space]
Consider the sub-Riemannian manifold given by $M=\R^3$, $\bD=\mathrm{span}\{X_1,X_2\}$,
$$
X_1=\partial_1,\quad X_2=\partial_2+\frac{x_1^2}{2}\partial_3,
$$
and the metric $dx_1^2+dx_2^2$. We choose $\varpi=dx_1\wedge dx_2\wedge dx_3$, that is, the canonical volume form on $\R^3$.

 The growth vector is equal to $(2,2,3)$ on the plane $N=\{ x_1=0\}$, and it is  $(2,3)$ elsewhere. As a consequence, $N$ is the set of singular points. At a regular point, $Q_{\mathrm{reg}}= 4$. Every singular point $p=(0,x_2,x_3)$ satisfies (A1) and we have $Q(p)=5$, $Q_N=4$, and $r_{\not N}=1$. Applying Corollaries~\ref{cor:dimHM} and~\ref{cor:infinite}, we obtain:
 $$
 \dim_H(M)=4, \qquad \hbox{and} \qquad \hh^4(B(p,\rho))<\infty \
 \hbox{ if $p$ regular,} \quad \hh^4(B(p,\rho))=\infty \ \hbox{ otherwise.}
 $$
  Thus small balls centered at singular points have infinite Hausdorff volume. This result can also be obtained by a direct computation based on the uniform Ball-Box Theorem, see~\cite{cimpa}.

  Note that the only family $(X_{I_{1}},X_{I_{2}},X_{I_{3}})$ such that $\sum_i |I_i| = Q_{\mathrm{reg}}$ is $(X_1, X_2, [X_1,X_2])$. The volume form of this family equals $x_1$ and it is of order 1 at every point of $N$. Thus every singular point satisfies assumptions (A1) and (A2) with $\sigma = 1$ ($\sigma = Q(p)-Q_{\mathrm{reg}}$ here).
\end{example}

\begin{example}
Consider the sub-Riemannian manifold given by $M=\R^4$, $\bD=\mathrm{span}\{X_1,X_2,X_3\}$, where
$$
X_1=\partial_1,\quad X_2=\partial_2+\frac{x_1^2}{2}\partial_4,\quad X_3=\partial_3+\frac{x_2^2}{2}\partial_4,
$$
and $g=dx_1^2+dx_2^2 +dx_4^2$. We choose $\varpi$ as the canonical volume form on $\R^4$.

At a regular point, $Q_{\mathrm{reg}}= 5$. The set of singular points is $N=\{ x_1=x_2=0\}$. Every singular point satisfies (A1) and we have $Q(p)=6$, $Q_N=4$, and $r_{\not N}=1$. Thus, by Corollary~\ref{cor:dimHM}, $\dim_H(M)=5$. However Corollary~\ref{cor:infinite} does not allow to conclude on the finiteness of the Hausdorff volume.

The only families such that $\sum_i |I_i| = Q_{\mathrm{reg}}$ are $(X_1, X_2, X_3, [X_1,X_2])$ and  $(X_1, X_2, X_3,[X_2,X_3])$. The volume form applied to these families is equal to $x_1$ and $x_2$ respectively,  and both of them are of order 1 at every point of $N$. Thus every singular point satisfies assumptions (A1) and (A2) with $\sigma = 1$ ($\sigma = Q(p)-Q_{\mathrm{reg}}$ here). Applying Theorem~\ref{th:finite}, we obtain:
 $$
 \dim_H(M)=5, \quad \hbox{and} \quad \hh^5(B(p,\rho))<\infty \ \hbox{ for any $p \in M$.}
 $$
\end{example}

\begin{example}
Let $M=\R^5$,  $\bD=\mathrm{span}\{X_1,X_2,X_3\}$,
$$
X_1= \partial_1, \quad X_2 = \partial_2 + x_1 \partial_3 + x_1^2 \partial_5, \quad X_3 = \partial_ 4 + x_1^k \partial_5,
$$
with $k>2$, and $g=dx_1^2+dx_2^2 +dx_3^2$. We choose $\varpi$ as the canonical volume form on $\R^5$.

The singular set is $N = \{x_1 = 0\}$.
A simple computation shows that every singular point $p$ satisfies (A1) and (A2), and $Q_{\mathrm{reg}}=7$, $Q(p)=8$, $Q_N=7$, $r_{\not N}=1$, and $\sigma = k-1$. Thus in this example $\sigma > Q(p)-Q_{\mathrm{reg}}$. Now Corollaries~\ref{cor:dimHM} and~\ref{cor:infinite} apply and we obtain
  $$
 \dim_H(M)=7, \qquad \hbox{and} \qquad \hh^7(B(p,\rho))<\infty \ \hbox{
   if $p$ regular,} \quad \hh^7(B(p,\rho))=\infty \ \hbox{ otherwise.}
 $$
\end{example}

\begin{example}
Let $M=\R^5$,  $\bD=\mathrm{span}\{X_1,X_2,X_3\}$,
$$
X_1= \partial_1, \quad X_2 = \partial_2 + x_1 \partial_3 + x_1^2 \partial_5, \quad X_3 = \partial_ 4 + (x_1^k + x_2^k) \partial_5,
$$
with $k>2$, and $g=dx_1^2+dx_2^2 +dx_3^2$. We choose $\varpi$ as the canonical volume form on $\R^5$.

The singular set is $N  = \{x_1 = x_2= 0\}$. Every singular point $p$ satisfies (A1) and (A2) and we have $Q_{\mathrm{reg}}=7$, $Q(0)=8$, $Q_N=6$, $r_{\not N}=1$, and $\sigma = k-1$. By Corollary~\ref{cor:dimHM} and Theorem~\ref{th:finite}, we obtain
$$
 \dim_H(M)=7, \qquad \hbox{and} \qquad \hh^7(B(p,\rho))<\infty \  \hbox{
   if $p$ regular,} \quad \hh^7(B(p,\rho))=\infty \ \hbox{ otherwise.}
 $$
Note that in this case we do not have $Q_{\mathrm{reg}}-Q_N < r_{\not N}$.
This shows that the criterion in Corollary~\ref{cor:infinite} does not provide a necessary condition for the Hausdorff volume to be infinite.
\end{example}

\begin{acknowledgement}
This work was supported by the European project AdG ERC ``GeMeThNES'', grant agreement number 246923  (see also {\tt gemethnes.sns.it}); by Digiteo grant {\it Congeo}; by the ANR project {\it GCM}, program ``Blanche'',
project number NT09\_504490; and by the Commission
of the European Communities under the 7th Framework Programme Marie
Curie Initial Training Network (FP7-PEOPLE-2010-ITN), project SADCO,
contract number 264735.
\end{acknowledgement}
%
%\section*{Appendix}
%\addcontentsline{toc}{section}{Appendix}
%
%

\bibliographystyle{abbrv}
\bibliography{biblio_hdim}

\begin{thebibliography}{10}

\bibitem{balu}
A.~Agrachev, D.~Barilari, and U.~Boscain.
\newblock On the {H}ausdorff volume in sub-{R}iemannian geometry.
\newblock {\em Calc. Var. Partial Differential Equations}, 43:355--388, 2012.

\bibitem{ambro}
L.~Ambrosio.
\newblock Fine properties of sets of finite perimeter in doubling metric
  measure spaces.
\newblock {\em Set-Valued Anal.}, 10(2-3):111--128, 2002.
\newblock Calculus of variations, nonsmooth analysis and related topics.

\bibitem{bellaiche}
A.~Bella{\"{\i}}che.
\newblock The tangent space in sub-{R}iemannian geometry.
\newblock In {\em Sub-{R}iemannian geometry}, volume 144 of {\em Progr. Math.},
  pages 1--78. Birkh\"auser, Basel, 1996.

\bibitem{fssc}
B.~Franchi, R.~Serapioni, and F.~Serra~Cassano.
\newblock On the structure of finite perimeter sets in step 2 {C}arnot groups.
\newblock {\em J. Geom. Anal.}, 13(3):421--466, 2003.

\bibitem{noi}
R.~Ghezzi and F.~Jean.
\newblock A new class of $(\hh^k,1)$-rectifiable subsets of metric spaces.
\newblock {\em Communications on Pure and Applied Analysis}, 12(2):881 -- 898,
  2013.

\bibitem{Goresky1988}
M.~Goresky and R.~MacPherson.
\newblock {\em Stratified {M}orse theory}, volume~14 of {\em Ergebnisse der
  Mathematik und ihrer Grenzgebiete (3) [Results in Mathematics and Related
  Areas (3)]}.
\newblock Springer-Verlag, Berlin, 1988.

\bibitem{gromov}
M.~Gromov.
\newblock {\em Structures m\'etriques pour les vari\'et\'es riemanniennes},
  volume~1 of {\em Textes Math\'ematiques [Mathematical Texts]}.
\newblock CEDIC, Paris, 1981.
\newblock Edited by J. Lafontaine and P. Pansu.

\bibitem{gromovcc}
M.~Gromov.
\newblock Carnot-{C}arath\'eodory spaces seen from within.
\newblock In {\em Sub-{R}iemannian geometry}, volume 144 of {\em Progr. Math.},
  pages 79--323. Birkh\"auser, Basel, 1996.

\bibitem{hermes}
H.~Hermes.
\newblock Nilpotent and high-order approximations of vector field systems.
\newblock {\em SIAM Rev.}, 33(2):238--264, 1991.

\bibitem{jean1}
F.~Jean.
\newblock Uniform estimation of sub-{R}iemannian balls.
\newblock {\em J. Dynam. Control Systems}, 7(4):473--500, 2001.

\bibitem{cimpa}
F.~{Jean}.
\newblock {Control of Nonholonomic Systems and Sub-Riemannian Geometry}.
\newblock {\em ArXiv e-prints}, 1209.4387, Sept. 2012.
\newblock Lectures given at the CIMPA School ``{G}\'eom\'etrie
  sous-riemannienne'', Beirut, Lebanon.

\bibitem{mitchell}
J.~Mitchell.
\newblock On {C}arnot-{C}arath\'eodory metrics.
\newblock {\em J. Differential Geom.}, 21(1):35--45, 1985.

\bibitem{montgomery}
R.~Montgomery.
\newblock {\em A tour of subriemannian geometries, their geodesics and
  applications}, volume~91 of {\em Mathematical Surveys and Monographs}.
\newblock American Mathematical Society, Providence, RI, 2002.

\end{thebibliography}

\end{document}